\documentclass[12pt,vatola]{article}

\usepackage{graphicx}

\def\c{\centerline}

\def\re#1{\par\hangindent\parindent\indent\llap{#1\enspace}\ignorespaces}

\def\no{\noindent}

\begin{document}

\c{\bf\large A New View of Combinatorial Maps}

\vskip 2mm

\c{\bf\large by Smarandache's Notion}

\vskip 5mm

\c{Linfan MAO}

\c{\scriptsize (Academy of Mathematics and Systems of Chinese
Academy of Sciences, Beijing 100080, P.R.China)}

\c{\scriptsize e-mail: maolinfan@163.com}

\vskip 6mm

\begin{minipage}{130mm}

\no{\bf Abstract}: {\small On a geometrical view, the conception
of map geometries is introduced, which is a nice model of the
Smarandache geometries, also new kind of and more general
intrinsic geometry of surfaces. Some open problems related
combinatorial maps with the Riemann geometry and Smarandache
geometries are presented.}

\no{\bf Key Words}: map, Smarandache geometry, model,
classification.

\no{\bf AMS(2000)}: 05C15, 20H15, 51D99, 51M05

\end{minipage}

\vskip 6mm

{\bf $1.$ What is a combinatorial map}

\vskip 3mm

A {\it graph} $\Gamma$ is a $2$-tuple $(V,E)$ consists of a finite
non-empty set $V$ of vertices together with a set $E$ of unordered
pairs of vertices, i.e., $E\subseteq V\times V$. Often denoted by
$V(\Gamma )$, $E(\Gamma )$ the vertex set and edge set of the
graph $\Gamma$([$9$]).

For example, the graph in the Fig.$1$ is the complete graph $K_4$
with vertex set $V=\{1,2,3,4\}$ and edge set $E=\{12,13,
14,23,24,34\}$.

A map is a connected topological graph cellularly embedded in a
surface. In 1973, Tutte gave an algebraic representation for an
embedding of a graph on locally orientable surface $([18])$, which
transfer a geometrical partition of a surface to a kind of
permutation in algebra as follows($[7][8]$).

A {\it{combinatorial map}} $M = ({\cal X}
_{\alpha,\beta},\cal{P})$ is defined to be a basic permutation
$\cal{P}$, i.e, for any $x\in {\cal X}_{\alpha,\beta}$, no integer
$k$ exists such that ${\cal{P}}^{k}x = \alpha x$, acting on ${\cal
X} _{\alpha,\beta}$, the disjoint union of {\it quadricells} $Kx$
of $x\in  X$ (the base set), where $K=\{1,\alpha,\beta,\alpha\beta
\}$ is the {\it Klein group}, satisfying the following two
conditions:\vskip 3mm

($i$) \ {\it  $\alpha{\cal{P}}={\cal{P}}^{-1}\alpha$;}\vskip 2mm

($ii$) \ {\it the group $\Psi_{J}=<\alpha,\beta,\cal{P}>$ is
transitive on ${\cal X}_{\alpha,\beta}$.}\vskip 2mm

For a given map $M=({\mathcal X}_{\alpha ,\beta},{\mathcal P})$,
it can be shown that $M^* = ({\mathcal X}_{\beta ,\alpha
},{\mathcal P}\alpha\beta)$ is also a map, call it the {\it dual}
of the map $M$. The vertices of $M$ are defined as the pairs of
conjugatcy orbits of ${\mathcal P}$ action on ${\mathcal
X}_{\alpha,\beta}$ by the condition $(Ci)$ and edges the orbits of
$K$ on ${\mathcal X}_{\alpha,\beta}$, for example,$\forall x\in
{\mathcal X}_{\alpha,\beta}$, $\{x,\alpha x,\beta x,\alpha\beta
x\}$ is an edge of the map $M$. Define the faces of $M$ to be the
vertices in the dual map $M^*$. Then the Euler characteristic
$\chi (M)$ of the map $M$ is

$$\chi (M)= \nu (M)-\varepsilon (M)+\phi (M)$$

\no{where,$\nu (M), \varepsilon (M), \phi (M)$ are the number of
vertices, edges and faces of the map $M$, respectively. For each
vertex of a map $M$, its valency is defined to be the length of
the orbits of ${\mathcal P}$ action on a quadricell incident with
$u$. }

For example, the graph $K_4$ on the tours with one face length $4$
and another $8$ , can be algebraic represented as follows:

A map $({\mathcal X}_{\alpha,\beta},\mathcal{P})$ with ${\mathcal
X}_{\alpha,\beta}= \{x,y,z,u,v,w,\alpha x,\alpha y, \alpha
z,\alpha u,\alpha v,\alpha w, \beta x,\beta y,\beta z,$ $\beta u,
 \beta v,\beta w,\alpha\beta x,
\alpha \beta y,\alpha \beta z,\alpha \beta u,\alpha \beta
v,\alpha\beta w \}$ and

\begin{eqnarray*}
{\mathcal P} &=& (x,y,z)(\alpha \beta x,u,w)(\alpha \beta z,\alpha
\beta u,v)
(\alpha \beta y,\alpha \beta v,\alpha \beta w)\\
&\times& (\alpha x,\alpha z,\alpha y)(\beta x,\alpha w,\alpha
u)(\beta z,\alpha v,\beta u)(\beta y,\beta w,\beta v)
\end{eqnarray*}

\no{The four vertices of this map are $ \{(x,y,z), (\alpha
x,\alpha z,\alpha y)\}$, $\{(\alpha \beta x,u,w),(\beta x,\alpha
w,\alpha u)\}$, $\{(\alpha \beta z,\alpha \beta u,v),(\beta
z,\alpha v,\beta u)\}$ and $\{(\alpha \beta y,\alpha \beta
v,\alpha \beta w),(\beta y,\beta w,\beta v)\}$ and six edges are
$\{e,\alpha e,\beta e,\alpha\beta e\}$, where, $e\in
\{x,y,z,u,v,w\}$. The Euler characteristic $\chi (M)$ is $\chi
(M)=4-6+2=0$.}

Geometrically, an embedding $M$ of a graph $\Gamma$ on a surface
is a map and has an algebraic representation. The graph $\Gamma$
is said the {\it underlying graph} of the map $M$ and denoted by
$\Gamma =\Gamma (M)$. For determining a given map $({\mathcal
X}_{\alpha,\beta},\mathcal{P})$ is orientable or not, the
following condition is needed.\vskip 3mm

 ($iii$) \ {\it If the group $\Psi_I=<\alpha\beta ,{\mathcal P}>$ is
transitive on ${\mathcal X}_{\alpha,\beta}$, then $M$ is
non-orientable. Otherwise, orientable.}\vskip 2mm

It can be shown that the number of orbits of the group
$\Psi_I=<\alpha\beta ,{\mathcal P}>$ in the Fig.$2$ action on
${\mathcal X}_{\alpha,\beta}= \{x,y,z,u,v,w,\alpha x,\alpha y,$ $
\alpha z,\alpha u,\alpha v,\alpha w, \beta x,\beta y,\beta z,\beta
u,$ $ \beta v,\beta w,\alpha\beta x, \alpha \beta y,\alpha \beta
z,\alpha \beta u,\alpha \beta v,\alpha\beta w \}$ is $2$. Whence,
it is an orientable map and the genus of the surface is $1$.
Therefore, the algebraic representation is correspondent with its
geometrical mean.

\vskip 8mm

{\bf $2.$ What are lost in combinatorial maps}

\vskip 5mm

As we known, mathematics is a powerful tool of sciences for its
unity and neatness, without any shade of mankind. On the other
hand, it is also a kind of aesthetics deep down in one's mind.
There is a famous proverb says that {\it only the beautiful things
can be handed down to today}, which is also true for the
mathematics.

Here, the term {\it unity} and {\it neatness} is relative and
local, also have various conditions. For acquiring the target,
many unimportant matters are abandoned in the process. Whether are
those matters in this time still unimportant in another time? It
is not true. That is why we need to think the question: {\it what
are lost in the classical mathematics?}

For example, a compact surface is topological equivalent to a
polygon with even number of edges by identifying each pairs of
edges along a given direction on it($[17]$). If label each pair of
edges by a letter $e,e\in {\mathcal E}$, a surface $S$ is also
identifying to a cyclic permutation such that each edge $e,e\in
{\mathcal E}$ just appears two times in $S$, one is $e$ and
another is $e^{-1}$. Let $a,b,c,\cdots$ denote the letters in
${\mathcal E}$ and $A,B,C,\cdots$ the sections of successive
letters  in linear order on a surface $S$ (or a string of letters
on $S$). Then, a surface can be represented as follows:

$$S=(\cdots , A,a,B,a^{-1},C,\cdots),$$

\no{where£¬$a\in {\mathcal E}$,$A,B,C$ denote a string of letters.
Define three elementary transformations as follows:}

\vskip 2mm $(O_1)\quad\quad (A,a,a^{-1},B)\Leftrightarrow (A,B);$

\vskip 2mm $(O_2)\quad\quad (i)\quad
(A,a,b,B,b^{-1},a^{-1})\Leftrightarrow (A,c,B,c^{-1}) ;$

\quad\quad\quad\quad $(ii)\quad (A,a,b,B,a,b)\Leftrightarrow
(A,c,B,c); $

\vskip 2mm $(O_3)\quad\quad (i)\quad
(A,a,B,C,a^{-1},D)\Leftrightarrow (B,a,A,D,a^{-1},C);$

\quad\quad\quad\quad $(ii)\quad (A,a,B,C,a,D)\Leftrightarrow
(B,a,A,C^{-1},a,D^{-1}).$ \vskip 2mm

\no If a surface $S_0$ can be obtained by the elementary
transformation $O_1$-$O_3$ from a surface $S$, it is said that $S$
{\it elementary equivalent} with $S_0$, denoted by
$S\sim_{El}S_0$.

We have known the following formula in $[8]$:

$(i)\quad
(A,a,B,b,C,a^{-1},D,b^{-1},E)\sim_{El}(A,D,C,B,E,a,b,a^{-1},b^{-1});$

$(ii)\quad (A,c,B,c)\sim_{El} (A,B^{-1},C,c,c);$

$(iii)\quad (A,c,c,a,b,a^{-1},b^{-1})\sim_{El} (A,c,c,a,a,b,b).$
\vskip 3mm

\no Then we can get the classification theorem of compact surface
as follows($[14]$):

{\it Any compact surface is homeomorphic to one of the following
standard surfaces:

($P_0$) The sphere: $aa^{-1}$;

($P_n$) The connected sum of $n,n\geq 1$, tori:

$$a_1b_1a_1^{-1}b_1^{-1}a_2b_2a_2^{-1}b_2^{-1}\cdots a_nb_na_n^{-1}b_n^{-1};$$

($Q_n$) The connected sum of $n,n\geq 1$, projective planes:

$$a_1a_1a_2a_2\cdots a_na_n.$$}

Generally, a combinatorial map is a kind of decomposition of a
surface. Notice that all the standard surfaces are just one face
map underlying an one vertex graph. By combinatorial view, a
combinatorial map is also a surface. But this assertion need more
clarifying. For example, see the tetrahedron graph $\Pi_4$ in the
$R^3$ and a map $K_4$ on the sphere.  Whether we can say it is the
sphere? Certainly NOT. Since any point $u$ on a sphere has a
neighborhood $N(u)$ homeomorphic to the open disc, therefore, all
angles incident with the point $1$ must all be $120^{\circ}$
degree on a sphere. But in $\Pi_4$, they are all $60^{\circ}$
degree. For making them topologically same, i.e., homeomorphism,
we must blow up the $\Pi_4$ to a sphere, as shown in the Fig.$3$.
Whence, for getting the classification theorem of compact
surfaces, we lose the {\it angle,area,
volume,distance,curvature},$\cdots$, etc, which are also lost in
the combinatorial maps.

{\bf Klein Erlanger Program} says that {\it any geometry is
finding invariant properties under the transformation group of
this geometry}. This is essentially the group action idea and
widely used in mathematics today. In the combinatorial maps, we
know the following problems are applications of the Klein Erlanger
Program:\vskip 2mm

 ($i$){\it to determine isomorphism maps or rooted maps;}

($ii$){\it to determine equivalent embeddings of a graph;}

($iii$){\it to determine an embedding whether exists;}

($iv$){\it to enumerate maps or rooted maps on a surface;}

($v$){\it to enumerate embeddings of a graph on a surface;}

($vi$) $\cdots$, etc. \vskip 2mm

All the problems are extensively investigated by researches in the
last century and papers related those problems are still appearing
frequently on the journals today. Then, {\it what are their
importance to classical mathematics?} and {\it what are their
contributions to science?} Those are the central topics of this
paper.

\vskip 5mm

{\bf $3.$ The Smarandache geometries}

\vskip 3mm

The {\it Smarandache geometries} is proposed by Smarandache in
1969 ($[16]$), which is a generalization of the classical
geometries, i.e., the Euclid, Lobachevshy-Bolyai-Gauss and
Riemannian geometries may be united altogether in the same space,
by some Smarandache geometries. These last geometries can be
either partially Euclidean and partially Non-Euclidean, or
Non-Euclidean. It seems that the Smarandache geometries are
connected with the {\it Relativity Theory} (because they include
the Riemann geometry in a subspace) and with the {\it Parallel
Universes} (because they combine separate spaces into one space)
too([$5$]). For a detail illustration, we need to consider the
classical geometries.

The axioms system of {\it Euclid geometry} are the following:

\vskip 3mm

(A1){\it there is a straight line between any two points.}

(A2){\it  a finite straight line can produce a infinite straight
line continuously.}

(A3){\it any point and a distance can describe a circle.}

(A4){\it all right angles are equal to one another.}

(A5){\it if a straight line falling on two straight lines make the
interior angles on the same side less than two right angles, then
the two straight lines, if produced indefinitely, meet on that
side on which are the angles less than the two right
angles.}\vskip 2mm

The axiom (A5) can be also replaced by: \vskip 3mm

(A5'){\it  given a line and a point exterior this line, there is
one line parallel to this line.}

\vskip 2mm

The Lobachevshy-Bolyai-Gauss geometry, also called hyperbolic
geometry, is a geometry with axioms $(A1)-(A4)$ and the following
axiom $(L5)$:\vskip 3mm

(L5) {\it there are infinitely many line parallels to a given line
passing through an exterior point.}\vskip 2mm

The Riemann geometry, also called elliptic geometry, is a geometry
with axioms $(A1)-(A4)$ and the following axiom $(R5)$:\vskip 3mm

{\it there is no parallel to a given line passing through an
exterior point.} \vskip 2mm

By the thought of Anti-Mathematics: not in a nihilistic way, but
in a positive one, i.e., banish the old concepts by some new ones:
their opposites, Smarandache introduced the {\it paradoxist
geometry, non-geometry, counter-projective geometry and
anti-geometry} in $[16]$ by contradicts the axioms $(A1)-(A5)$ in
Euclid geometry, generalize the classical geometries.

\vskip 4mm

{\bf Paradoxist geometry}

\vskip 3mm

In this geometry, its axioms are $(A1)-(A4)$ and  with one of the
following as the axiom $(P5)$: \vskip 3mm

($i$){\it there are at least a straight line and a point exterior
to it in this space for which any line that passes through the
point intersect the initial line.}

($ii$){\it there are at least a straight line and a point exterior
to it in this space for which only one line passes through the
point and does not intersect the initial line.}

($iii$){\it there are at least a straight line and a point
exterior to it in this space for which only a finite number of
lines $l_1,l_2,\cdots , l_k, k\geq 2$ pass through the point and
do not intersect the initial line.}

($iv$){\it there are at least a straight line and a point exterior
to it in this space for which an infinite number of lines pass
through the point (but not all of them) and do not intersect the
initial line.}

($v$){\it there are at least a straight line and a point exterior
to it in this space for which any line that passes through the
point and does not intersect the initial line.}\vskip 2mm

\vskip 4mm

{\bf Non-Geometry}

\vskip 3mm

The non-geometry is a geometry by denial some axioms of
$(A1)-(A5)$, such as:

\vskip 3mm

($A1^-$){\it It is not always possible to draw a line from an
arbitrary point to another arbitrary point.}

($A2^-$){\it It is not always possible to extend by continuity a
finite line to an infinite line.}

($A3^-$){\it It is not always possible to draw a circle from an
arbitrary point and of an arbitrary interval.}

($A4^-$){\it not all the right angles are congruent.}

($A5^-$){\it if a line, cutting two other lines, forms the
interior angles of the same side of it strictly less than two
right angle, then not always the two lines extended towards
infinite cut each other in the side where the angles are strictly
less than two right angle.}

\vskip 4mm

{\bf Counter-Projective geometry}

\vskip 3mm

Denoted by $P$ the point set, $L$ the line set and $R$ a relation
included in $P\times L$. A counter-projective geometry is a
geometry with the following counter-axioms:\vskip 2mm

($C1$){\it There exist: either at least two lines, or  no line,
that contains two given distinct points.}

($C2$){\it Let $p_1,p_2,p_3$ be three non-collinear points, and
$q_1,q_2$ two distinct points. Suppose that $\{p_1.q_1,p_3\}$ and
$\{p_2,q_2,p_3\}$ are collinear triples. Then the line containing
$p_1,p_2$ and the line containing $q_1,q_2$ do not intersect.}

($C3$){\it Every line contains at most two distinct points. }

\vskip 4mm

{\bf Anti-Geometry}

\vskip 3mm

A geometry by denial some axioms of the Hilbert's $21$ axioms of
Euclidean geometry. As shown in $[5]$, there are at least
$2^{21}-1$ anti-geometries.

The Smarandache geometries are defined as follows.

\vskip 3mm

\no{\bf Definition $3.1$} {\it An axiom is said Smarandachely
denied if the axiom behaves in at least two different ways within
the same space, i.e., validated and invalided, or only invalided
but in multiple distinct ways.

A Smarandache geometry is a geometry which has at least one
Smarandachely denied axiom($1969$).}

\vskip 2mm

A nice model for the Smarandache geometries, called $s$-manifolds,
is found by Iseri in $[3][4]$, which is defined as follows:

{\it An $s$-manifold is any collection ${\mathcal C}(T,n)$ of
these equilateral triangular disks $T_i, 1\leq i\leq n$ satisfying
the following conditions:}

$(i)$ {\it Each edge $e$ is the identification of at most two
edges $e_i,e_j$ in two distinct triangular disks $T_i,T_j, 1\leq
i,j\leq n$ and $i\not= j$;}

$(ii)$ {\it Each vertex $v$ is the identification of one vertex in
each of five, six or seven distinct triangular disks.}

The vertices are classified by the number of the disks around
them. A vertex around five, six or seven triangular disks is
called an {\it elliptic vertex}, a {\it Euclid vertex} or a {\it
hyperbolic vertex}, respectively.

An $s$-manifold is called closed if each edge is shared by exactly
two triangular disks. An elementary classification for closed
$s$-manifolds by triangulation are made in the reference $[11]$.
The closed $s$-manifolds are classified into $7$ classes in
$[11]$, as follows:\vskip 3mm

{\bf Classical Type}:\vskip 2mm

$(1)$  $\Delta_1=\{5-regular \ triangular \ maps\}$ ({\it
elliptic});

$(2)$ $\Delta_2=\{6-regular  \ triangular \ maps\}$({\it
euclidean});

$(3)$ $\Delta_3=\{7-regular  \ triangular \ maps\}$({\it
hyperbolic}).\newpage

{\bf Smarandache Type}:\vskip 2mm

$(4)$ $\Delta_4=\{triangular \ maps \ with \ vertex \ valency \ 5
\ and \ 6\}$ ({\it euclid-elliptic});

$(5)$ $\Delta_5=\{triangular \ maps \ with \ vertex \ valency \ 5
\ and \ 7\}$ ({\it elliptic-hyperbolic});

$(6)$ $\Delta_6=\{triangular \ maps \ with \ vertex \ valency \ 6
\ and \ 7\}$ ({\it euclid-hyperbolic});

$(7)$ $\Delta_7=\{triangular \ maps \ with \ vertex \ valency \ 5
, 6 \ and \ 7\}$ ({\it mixed}).\vskip 2mm

It is proved in $[11]$ that $|\Delta_1|=2$, $|\Delta_5|\geq 2$ and
$|\Delta_i|, i=2,3,4,6,7$ are infinite. Isier proposed a question
in $[3]$: {\it Do the other closed $2$-manifolds correspond to
$s$-manifolds with only hyperbolic vertices?}. Since there are
infinite Hurwitz maps, i.e., $|\Delta_3|$ is infinite, the answer
is affirmative.

\vskip 8mm

{\bf $4.$ The map geometries}

\vskip 5mm

Combinatorial maps can be used to construct new geometries, which
are nice models for the Smarandache geometries, also a
generalization of Isier's model and Poincar\'{e}'s model for
hyperbolic geometry.

\vskip 4mm

{\bf $4.1$ \ Map geometries without boundary}

\vskip 3mm

For a given map on a surface, the map geometries without boundary
are defined as follows.

\vskip 3mm

\no{\bf Definition $4.1$} {\it For a combinatorial map $M$ with
each vertex valency$\geq 3$, associates a real number $\mu (u), 0
\ < \mu (u) \ < \ \pi$, to each vertex $u, u\in V(M)$. Call
$(M,\mu)$ a map geometry with out boundary, $\mu (u)$ the angle
factor of the vertex $u$ and to be orientablle or non-orientable
if $M$ is orientable or not.}

\vskip 2mm

The realization of each vertex $u, u\in V(M)$ in $R^3$ space is
shown in the Fig.$1$ for each case of $\rho (u)\mu (u) \ > \
2\pi$, $=2\pi$ or $ \ < \ 2\pi$.

\includegraphics[bb=20 20 100 140]{g4.eps}

\vskip 2mm

\c{$\rho (u)\mu (u) \ < \ 2\pi$\hskip 15mm $\rho (u)\mu (u)
=2\pi$\hskip 20mm $\rho (u)\mu (u) \ > \ 2\pi$}\vskip 2mm

\c{\bf Fig.$1$}\vskip 2mm

As pointed out in the Section $2$, this kind of realization is not
a surface, but it is homeomorphic to a surface. We classify points
in a map geometry $(M,\mu)$ with out boundary as follows.

\vskip 3mm

\no{\bf Definition $4.2$} {\it A point $u$ in a map geometry
$(M,\mu)$ is called elliptic, euclidean or hyperbolic if $\rho
(u)\mu (u) \ < \ 2\pi$, $\rho (u)\mu (u) =2\pi$ or $\rho (u)\mu
(u) \ > \ 2\pi$.}

Then we have the following results.

\vskip 3mm

\no{\bf Proposition $4.1$} {\it Let $M$  be a map with $\forall
u\in V(M), \rho (u)\geq 3$. Then for $\forall u\in V(M)$, there is
a map geometries $(M,\mu)$ without boundary such that $u$ is
elliptic, euclidean or hyperbolic in this geometry.}

\vskip 2mm

{\it Proof} \ Since $\rho (u)\geq 3$, we can choose the angle
factor $\mu (u)$ such that $\mu (u)\rho (u) \ < \ 2\pi$, $\mu
(u)\rho (u)= 2\pi$ or $\mu (u)\rho (u) \ > \ 2\pi$. Notice that

$$0 \ < \frac{2\pi}{\rho (u)} \ < \ \pi.$$

\no Whence, we can also choose $\mu (u)$ satisfying that $ 0 \ <
\mu (u) \ < \pi \quad\quad \natural$

\vskip 3mm

\no{\bf Proposition $4.2$} {\it Let $M$  be a map of order$\geq 3$
and $\forall u\in V(M), \rho (u)\geq 3$. Then there exists a map
geometry $(M,\mu)$ with out boundary, in which all points are one
of the elliptic vertices, euclidean vertices and hyperbolic
vertices or their mixed.}

\vskip 2mm

{\it Proof} \ According to the Proposition $4.1$,  we can choose
an angle factor $\mu$ such that a vertex $u, u\in V(M)$ to be
elliptic, or euclidean, or hyperbolic. Since $|V(M)|\geq 3$, we
can also choose the angle factor $\mu$ such that any two vertices
$v,w\in V(M)\backslash \{u\}$ to be elliptic, or euclidean, or
hyperbolic according to our wish. Then the map geometry $(M,\mu)$
makes the assertion hold.\quad\quad $\natural$

A {\it geodesic} in a manifold is a curve as straight as possible.
Similarly, in a map geometry, its $m$-lines and $m$-points are
defined as follows.

\vskip 3mm

\no{\bf Definition $4.3$} {\it Let $(M,\mu)$ be a map geometry
without boundary. An $m$-line in $(M,\mu)$ is a curve with a
constant curvature and points in it are called $m$-points.}

If an $m$-line pass through an elliptic point or a hyperbolic
point $u$, it must has the angle $\frac{\mu (u)\rho (u)}{2}$ with
the entering line, not $180^{\circ}$, which are explained in the
Fig.$2$.

\includegraphics[bb=60 5 100 150]{g6.eps}

\vskip 2mm

{\hskip 10mm ${\rm a}=\frac{\mu (u)\rho (u)}{2} \ <  \ \pi$ \hskip
60mm ${\rm a}=\frac{\mu (u)\rho (u)}{2} \ >  \ \pi$}\vskip 2mm

\c{\bf Fig.$2$}

\vskip 3mm

The following proposition asserts that all map geometries without
boundary are Smarandache geometries.

\vskip 4mm

\no{\bf Proposition $4.3$} {\it For a map $M$ on a locally
orientable surface with order$\geq 3$ and vertex valency$\geq 3$,
there is an angle factor $\mu$ such that $(M,\mu)$ is a
Smarandache geometry by denial the axiom (A5) with the axioms
(A5),(L5) and (R5). }

\vskip 3mm

{\it Proof} \ According to the Proposition $4.1$, we know that
there exist an angle factor $\mu$ such that there are elliptic
vertices, euclidean vertices and hyperbolic vertices in $(M,\mu)$
simultaneously. The proof is divided into three cases.

\vskip 4mm

{\bf Case $1.$ \ $M$ is a planar map}

\vskip 3mm

Notice that for a given line $L$ not pass through the vertices in
the map $M$ and a point $u$ on its left side in $(M,\mu)$, if $u$
is an euclidean point, then there is one and only one line passes
through $u$ not intersect with $L$, and  if $u$ is an elliptic
point, then there are infinite lines pass through $u$ not
intersect with $L$, but if $u$ is a hyperbolic point, then each
line passes through $u$ will intersect with $L$. Therefore,
$(M,\mu)$ is a Smarandache geometry by denial the axiom (A5) with
the axioms (A5), (L5) and (R5).

\vskip 4mm

{\bf Case $2.$ \ $M$ is an orientable map}

\vskip 3mm

According to the classification theorem of compact surfaces, We
only need to prove this result for the torus. Notice that on the
torus, an $m$-line has the following properties ([$15$]):

\vskip 3mm

{\it If the slope $\varsigma$ of $m$-line $L$ is a rational
number, then $L$ is a closed line on the torus. Otherwise, $L$ is
infinite, and moreover $L$ passes arbitrarily close to every point
of the torus.}

\vskip 3mm

Whence, if $L_1$ is an $m$-line on the torus, not passes through
an elliptic or hyperbolic point, then for any point $u$ exterior
$L_1$, we know that if $u$ is an euclidean point, then there is
only one $m$-line passes through $u$ not intersect with $L_1$, and
if $u$ is elliptic or hyperbolic, then any $m$-line passes through
$u$ will intersect with $L_1$.

Now let $L_2$ be an $m$-line passes through an elliptic or
hyperbolic point, such as the $m$-line in the Fig.$3$ and $v$ an
euclidean point.

\includegraphics[bb=60 5 100 180]{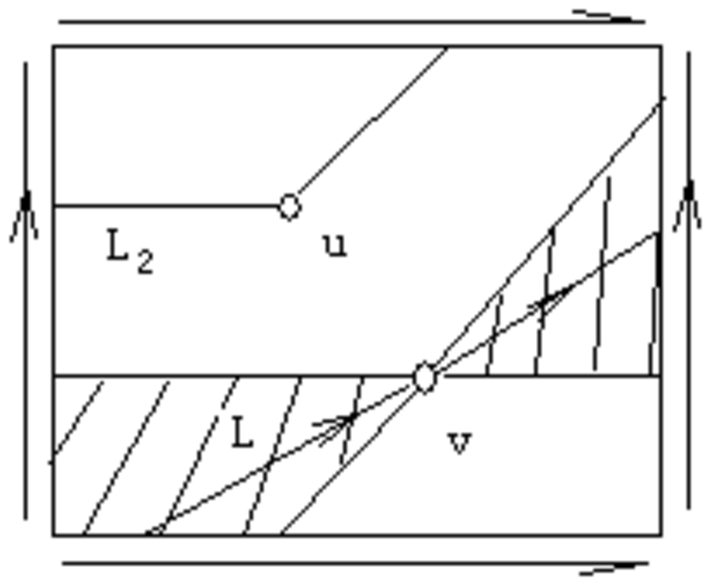}

\vskip 3mm

\c{\bf Fig.$3$}

\no Then any $m$-line $L$ in the shade filed passes through the
point $v$ will not intersect with $L_2$. Therefore, $(M,\mu)$ is a
Smarandache geometry by denial the axiom (A5) with the axioms
(A5),(L5) and (R5).

\vskip 3mm

{\bf Case $3.$ \ $M$ is a non-orientable map}

\vskip 2mm

Similar to the Case $2$, by the classification theorem of the
compact surfaces, we only need to prove this result for the
projective plane. Now let the $m$-line passes through the center
in the circle. Then if $u$ is an euclidean point, there is only
one $m$-line passes through $u$, see (a) in the Fig.$4$. If $v$ is
an elliptic point and there is an $m$-line passes through it and
intersect with $L$, see (b) in the Fig.$4$, assume the point $1$
is a point such that the $m$-line $1v$ passes through $0$, then
any $m$-line in the shade of (b) passes through the point $v$ will
intersect with $L$.

\includegraphics[bb=20 5 100 200]{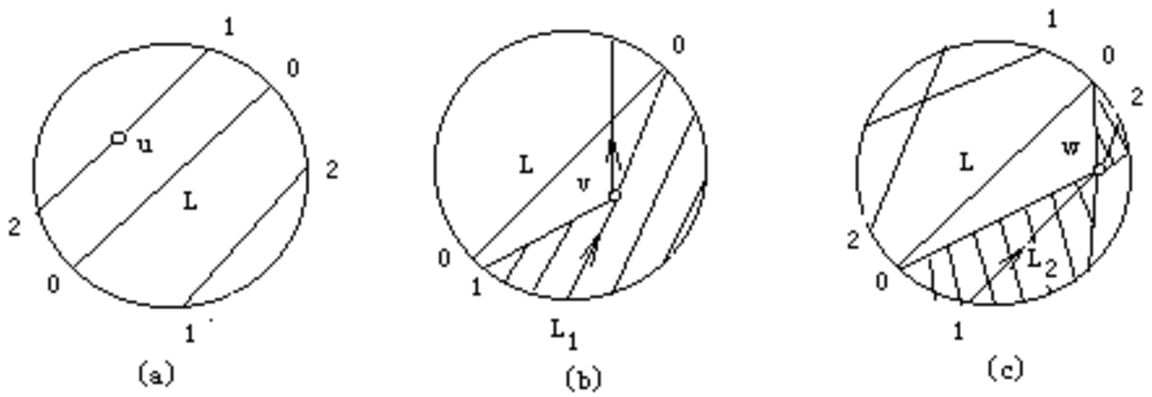}

\vskip 3mm

\c{\bf Fig.$4$}\vskip 2mm

If $w$ is a point and there is an $m$-line passes through it and
does not intersect with $L$, see (c) in the Fig.$4$, then any
$m$-line in the shade of (c) passes through the point $w$ will not
intersect with $L$. Since the position of the vertices of the map
$M$ on the projective plane can be choose as our wish, the proof
is complete. \quad\quad $\natural$.

\vskip 3mm

{\bf $4.2$ \ Map geometries with boundary}

\vskip 2mm

The Poincar\'{e}'s model for hyperbolic geometry hints us to
introduce the map geometries with boundary, which is defined as
follows.

\vskip 3mm

\no{\bf Definition $4.4$} {\it For a map geometry $(M,\mu )$
without boundary and faces $f_1,f_2,\cdots ,f_l\in F(M), 1\leq
l\leq \phi (M)-1$, if $(M,\mu)\setminus\{f_1,f_2,\cdots ,f_l\}$ is
connected, then call
$(M,\mu)^{-l}=(M,\mu)\setminus\{f_1,f_2,\cdots ,f_l\}$ a map
geometry with boundary $f_1,f_2,\cdots ,f_l$ and orientable or not
if $(M,\mu)$ is orientable or not.

A connected curve with constant curvature in $(M,\mu)^{-l}$ is
called an $m^-$-line and points $m^-$-points.}

The map geometries with boundary also are Smarandache geometries,
which is convince by the following result.

\vskip 3mm

\no{\bf Proposition $4.4$} {\it For a map $M$ on a locally
orientable surface with order$\geq 3$, vertex valency$\geq 3$ and
a face $f\in F(M)$, there is an angle factor $\mu$ such that
$(M,\mu)^{-1}$ is a Smarandache geometry by denial the axiom (A5)
with the axioms (A5),(L5) and (R5). }

\vskip 2mm

{\it Proof} \ Similar to the proof of the Proposition $4.3$,
consider the map $M$ being a planar map, an orientable map on a
torus or a non-orientable map on a projective plane, respectively.
We get the assertion. \quad\quad $\natural$

Notice that for an one face map geometry $(M,\mu)^{-1}$ with
boundary, if we choose all points being euclidean, then
$(M,\mu)^{-1}$ is just the Poincar\'{e}'s model for hyperbolic
geometry.

\vskip 4mm

{\bf $4.3$ \ Classification of map geometries}

\vskip 3mm

For the classification of the map geometries, we introduce the
following definition.

\vskip 3mm

\no{\bf Definition $4.5$} {\it Two map geometries $(M_1,\mu_1)$
and $(M_2,\mu_2)$ or $(M_1,\mu_1)^{-l}$ and $(M_2,\mu_2)^{-l}$ are
equivalent if there is a bijection $\theta :M_1\rightarrow M_2$
such that for $\forall u\in V(M)$, $\theta (u)$ is euclidean,
elliptic or hyperbolic iff $u$ is euclidean, elliptic or
hyperbolic.}

\vskip 2mm

The relation of the numbers of unrooted maps with the map
geometries is the following.

\vskip 4mm

\no{\bf Proposition $4.5$} {\it If $\mathcal M$ is a set of
non-isomorphisc maps with order $n$ and $m$ faces, then the number
of map geometries without boundary is $3^n|\mathcal M|$ and the
number of map geometries with one face being its boundary is
$3^nm|\mathcal M|$.}

\vskip 3mm

{\it Proof} \ By the definition, for a map $M\in {\mathcal M}$,
there are $3^n$ map geometries without boundary and $3^nm$ map
geometries with one face being its boundary by the Proposition
$4.3$. Whence, we get $3^n|\mathcal M|$ map geometries without
boundary and $3^nm|\mathcal M|$ map geometries with one face being
its boundary from $\mathcal M$.\quad\quad $\natural$.

We have the following enumeration result for the non-equivalent
map geometries without boundary.

\vskip 3mm

\no{\bf Proposition $4.6$} {\it The numbers $n^O(\Gamma ,g)$,
$n^N(\Gamma ,g)$ of non-equivalent orientable, non-orientable map
geometries without boundary underlying a simple graph $\Gamma$ by
denial the axiom (A5) by (A5), (L5) or (R5) are  }

$$n^O(\Gamma ,g)=\frac{3^{|\Gamma|}\prod\limits_{v\in V(\Gamma)}(\rho (v)-1)!}{2|{\rm Aut}\Gamma |},$$

\no{\it and}

$$n^N(\Gamma ,g)=\frac{(2^{\beta (\Gamma )}-1)3^{|\Gamma|}\prod\limits_{v\in V(\Gamma)}(\rho (v)-1)!}{2|{\rm Aut}\Gamma |},$$

\no{\it where $\beta (\Gamma)=\varepsilon (\Gamma)-\nu (\Gamma)
+1$ is the  Betti number of the graph $\Gamma$.}

\vskip 2mm

{\it Proof} Denote by ${\mathcal M}(\Gamma)$  the set of
non-isomorphic maps underlying the graph $\Gamma$ on locally
orientable surfaces and by ${\mathcal E}(\Gamma)$ the set of
embeddings of the graph $\Gamma$ on the locally orientable
surfaces. For a map $M, M\in {\mathcal M}(\Gamma)$, there are
$\frac{3^{|M|}}{|{\rm Aut}M|}$ different map geometries without
boundary by choosing the angle factor $\mu$ on a vertex $u$ such
that $u$ is euclidean, elliptic or hyperbolic. From permutation
groups, we know that

$$
|{\rm Aut}\Gamma\times <\alpha>|=|({\rm Aut}\Gamma)_M||M^{{\rm
Aut}\Gamma\times <\alpha>}|=|{\rm Aut}M||M^{{\rm Aut}\Gamma\times
<\alpha>}|.
$$

\no Therefore, we get that

\begin{eqnarray*}
n^O(\Gamma ,g)& = & \sum\limits_{M\in {\mathcal
M}(\Gamma)}\frac{3^{|M|}}{|{\rm Aut}M|}\\
&=& \frac{3^{|\Gamma |}}{|\rm Aut\Gamma\times
<\alpha>|}\sum\limits_{M\in {\mathcal M}(\Gamma)}\frac{|\rm
Aut\Gamma\times <\alpha>|}{|{\rm Aut}M|}\\
&=& \frac{3^{|\Gamma |}}{|\rm Aut\Gamma\times
<\alpha>|}\sum\limits_{M\in {\mathcal M}(\Gamma)}|M^{{\rm
Aut}\Gamma\times <\alpha>}|\\
&=& \frac{3^{|\Gamma |}}{|\rm Aut\Gamma\times <\alpha>|}|{\mathcal
E}^O(\Gamma)|\\
&=& \frac{ 3^{|\Gamma|}\prod\limits_{v\in V(\Gamma)}(\rho
(v)-1)!}{2|{\rm Aut}\Gamma |}.
\end{eqnarray*}

\no Similarly, we get that

\begin{eqnarray*}
n^N(\Gamma ,g)& = & \frac{3^{|\Gamma |}}{|\rm Aut\Gamma\times
<\alpha>|}|{\mathcal
E}^N(\Gamma)|\\
&=& \frac{ (2^{\beta (\Gamma )}-1)3^{|\Gamma|}\prod\limits_{v\in
V(\Gamma)}(\rho (v)-1)!}{2|{\rm Aut}\Gamma |}.
\end{eqnarray*}

\no This completes the proof. \quad\quad $\natural$

For the classification of map geometries with boundary, we have
the following result.

\vskip 3mm

\no{\bf Proposition $4.7$} {\it The numbers $n^O(\Gamma ,-g)$,
$n^N(\Gamma ,-g)$ of non-equivalent orientable, non-orientable map
geometries with one face being its boundary and underlying a
simple graph $\Gamma$ by denial the axiom (A5) by (A5), (L5) or
(R5) are}

$$n^O(\Gamma ,-g)=\frac{3^{|\Gamma|}}{2|{\rm Aut}\Gamma|}[(\beta (\Gamma)+1)
\prod\limits_{v\in V(\Gamma)}(\rho (v)-1)!-\frac{2d(g[\Gamma
](x))}{dx}|_{x=1}]$$

\no{\it and}

$$n^N(\Gamma ,-g)=\frac{(2^{\beta (\Gamma)}-1)3^{|\Gamma|}}{2|{\rm Aut}\Gamma|}[(\beta (\Gamma)+1)
\prod\limits_{v\in V(\Gamma)}(\rho (v)-1)!-\frac{2d(g[\Gamma
](x))}{dx}|_{x=1}],$$

\no{\it where $g[\Gamma ](x)$ is the genus polynomial of the graph
$\Gamma$ ( see {\rm [$12$]}), i.e., $g[\Gamma
](x)=\sum\limits_{k=\gamma (\Gamma)}^{\gamma_m(\Gamma)}g_k[\Gamma
]x^k$ with $g_k[\Gamma ]$ being the number of embeddings of
$\Gamma$ on the orientable surface of genus $k$.}

\vskip 2mm

{\it Proof} Notice that $\nu (M)-\varepsilon (M)+\phi (M)=2-2g(M)$
for an orientable map $M$ by the Euler characteristic. Similar to
the proof of the Proposition $4.6$ with the notation ${\mathcal
M}(\Gamma)$, by the Proposition $4.5$ we know that

\begin{eqnarray*}
n^O(\Gamma ,-g)&=& \sum\limits_{M\in{\mathcal
M}(\Gamma)}\frac{\phi (M)3^{|M|}}{|{\rm Aut}M|}\\
&=&\sum\limits_{M\in{\mathcal M}(\Gamma)}\frac{(2+\varepsilon
(\Gamma)-\nu (\Gamma)-2g(M))3^{|M|}}{|{\rm Aut}M|}\\
&=&\sum\limits_{M\in{\mathcal M}(\Gamma)}\frac{(2+\varepsilon
(\Gamma)-\nu (\Gamma))3^{|M|}}{|{\rm
Aut}M|}-\sum\limits_{M\in{\mathcal
M}(\Gamma)}\frac{2g(M)3^{|M|}}{|{\rm Aut}M|}\\
&=&\frac{(2+\varepsilon (\Gamma)-\nu (\Gamma))3^{|M|}}{|{\rm
Aut}\Gamma\times <\alpha>|}\sum\limits_{M\in{\mathcal
M}(\Gamma)}\frac{|{\rm Aut}\Gamma\times <\alpha>|}{|{\rm
Aut}M|}\\
&-&\frac{2\times 3^{|\Gamma |}}{|{\rm Aut}\Gamma\times
<\alpha>|}\sum\limits_{M\in{\mathcal M}(\Gamma)}\frac{g(M)|{\rm
Aut}\Gamma\times
<\alpha>|}{|{\rm Aut}M|}\\
&=& \frac{(\beta (\Gamma)+1)3^{|M|}}{|{\rm Aut}\Gamma\times
<\alpha>|}\sum\limits_{M\in{\mathcal M}}(\Gamma)|M^{{\rm
Aut}\Gamma\times <\alpha>}|\\
&-&\frac{3^{|\Gamma |}}{|{\rm Aut}\Gamma
|}\sum\limits_{M\in{\mathcal M}(\Gamma)}g(M)|M^{{\rm
Aut}\Gamma\times <\alpha>}|\\
&=&\frac{(\beta (\Gamma)+1)3^{|\Gamma |}}{2|{\rm Aut}\Gamma
|}\prod\limits_{v\in V(\Gamma)}(\rho (v)-1)! -\frac{3^{|\Gamma
|}}{|{\rm Aut}\Gamma |}\sum\limits_{k=\gamma
(\Gamma)}^{\gamma_m(\Gamma)}kg_k[\Gamma ]\\
&=&\frac{3^{|\Gamma |}}{2|{\rm Aut}\Gamma |}[(\beta
(\Gamma)+1)\prod\limits_{v\in V(\Gamma)}(\rho (v)-1)!-
\frac{2d(g[\Gamma ](x))}{dx}|_{x=1}].
\end{eqnarray*}

Notice that $n^L(\Gamma ,-g)=n^O(\Gamma ,-g)+n^N(\Gamma ,-g)$ and
the number of re-embeddings of an orientable map $M$ on surfaces
is $2^{\beta (M)}$ (see also [$13$]). We have that

\begin{eqnarray*}
n^L(\Gamma ,-g)&=&\sum\limits_{M\in {\mathcal
M}(\Gamma)}\frac{2^{\beta (M)}\times 3^{|M|}\phi (M)}{|{\rm
Aut}M|}\\
&=& 2^{\beta (M)}n^O(\Gamma ,-g).
\end{eqnarray*}

\no Whence, we get that

\begin{eqnarray*}
n^N(\Gamma ,-g)&=& (2^{\beta (M)}-1)n^O(\Gamma ,-g)\\
&=& \frac{(2^{\beta (M)}-1)3^{|\Gamma |}}{2|{\rm Aut}\Gamma
|}[(\beta (\Gamma)+1)\prod\limits_{v\in V(\Gamma)}(\rho (v)-1)!-
\frac{2d(g[\Gamma ](x))}{dx}|_{x=1}].
\end{eqnarray*}

\no This completes the proof.\quad\quad $\natural$

\vskip 4mm

{\bf $4.4$ \ Polygons in a map geometry}

\vskip 3mm

A $k$-polygon in a map geometry is a $k$-polygon with each line
segment being $m$-lines or $m^-$-lines. For the sum of the
internal angles in a $k$-polygon, we have the following result.

\vskip 3mm

\no{\bf Proposition $4.8$} {\it Let $P$ be a $k$-polygon in a map
geometry with each line segment passes through at most one
elliptic or hyperbolic point. If $H$ is the set of elliptic points
and hyperbolic points on the line segment of $P$, then the sum of
the internal angles in $P$ is}

$$(k+|H|-2)\pi-\frac{1}{2}\sum\limits_{u\in H}\rho (u)\mu (u).$$

\vskip 2mm

{\it Proof} \ Denote by $U, V$ the sets of elliptic points and
hyperbolic points in $H$ and $|U|=p, |V|=q$.  If an $m$-line
segment passes through an elliptic point $u$, add a straight line
segment in the plane as the Fig.$6$(1). Then we get that

$${\rm angle} \  {\rm a}= {\rm angle} 1 + {\rm angle} 2=\pi - \frac{\rho (u)\mu (u)}{2}.$$

If an $m$-line passes through an hyperbolic point $v$, also add a
straight line segment in the plane as the Fig.$6$(2). Then we get
that

$${\rm angle} \ b= {\rm angle} 3 + {\rm angle} 4=\frac{\rho (v)\mu (v)}{2}-\pi.$$

\vskip 3mm

\includegraphics[bb=60 5 100 180]{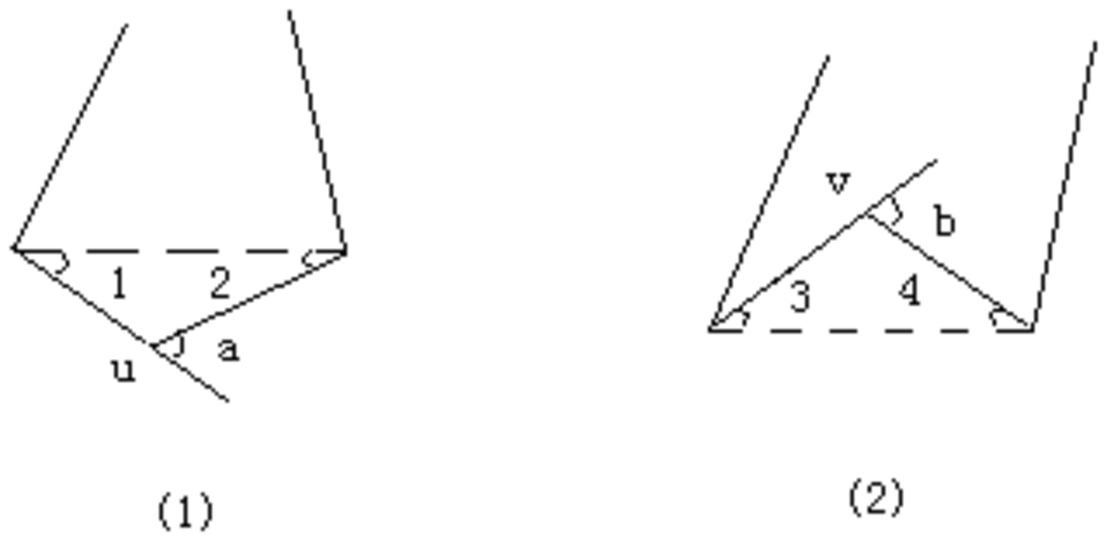}

\vskip 2mm

\c{\bf Fig.$5$}\vskip 2mm

Since the sum of the internal angles of a $k$-polygon in the plane
is $(k-2)\pi$, we know that the sum of the internal angles in $P$
is

\begin{eqnarray*}
&(k& -2)\pi+\sum\limits_{u\in U}(\pi-\frac{\rho (u)\mu
(u)}{2})-\sum\limits_{v\in V}(\frac{\rho (u)\mu (u)}{2}-\pi)\\
&=& (k+p+q-2)\pi-\frac{1}{2}\sum\limits_{u\in H}\rho (u)\mu
(u)\\
&=& (k+|H|-2)\pi-\frac{1}{2}\sum\limits_{u\in H}\rho (u)\mu (u).
\end{eqnarray*}

\no This completes the proof. \quad\quad $\natural$

As a corollary, we get the sum of the internal angles of a
triangle in a map geometry as follows, which is consistent with
the classical results.

\vskip 3mm

\no{\bf Corollary $4.1$} {\it Let $\triangle$ be a triangle in a
map geometry. Then

($i$) if $\triangle$ is euclidean, then then the sum of its
internal angles is equal to $\pi$;

($ii$) if $\triangle$ is elliptic, then the sum of its internal
angles is less than $\pi$;

($iii$) if $\triangle$ is hyperbolic, then the sum of its internal
angles is more than $\pi$.
 }

\vskip 8mm

{\bf $5.$ Open problems for applying maps to classical geometries}

\vskip 5mm

Here is a collection of open problems concerned combinatorial maps
with the Riemann geometry and Smarandache geometries. Although
they are called open problems, in fact, any solution for one of
these problems needs to establish a new mathematical system first.

\vskip 4mm

{\bf $5.1$ The uniformization theorem for simple connected Riemann
surfaces}

\vskip 3mm

The {\it uniformization theorem} for simple connected Riemann
surfaces is one of those beautiful results in the Riemann surface
theory, which is stated as follows([$2$]).

\vskip 3mm

{\it If ${\mathcal S}$ is a simple connected Riemann surface, then
${\mathcal S}$ is conformally equivalent to one and only one of
the following three:

$(a)$\quad ${\mathcal C}\bigcup {\infty}$;

$(b)$\quad ${\mathcal C}$;

$(c)$\quad $\triangle =\{z\in {\mathcal C}| |z| \ < 1\}.$}\vskip
2mm

\no{We have proved in  $[11]$ that any automorphism of a map is
conformal. Therefore, we can also introduced the conformal mapping
between maps. Then, {\it how can we define the conformal
equivalence for maps enabling us to get the uniformization theorem
of maps? What is the correspondence class maps with the three type
$(a)-(c)$ Riemann surfaces?}}

\vskip 4mm

{\bf $5.2$ Combinatorial construction of an algebraic curve of
genus}

\vskip 3mm

A {\it complex plane algebraic curve ${\mathcal C}_l$} is a
homogeneous equation $f(x,y,z)=0$ in $P_2{\mathcal
C}=(C^2\setminus (0,0,0))/\sim$, where $f(x,y,z)$ is a polynomial
in $x,y$ and $z$ with coefficients in ${\mathcal C}$. The degree
of $f(x,y,z)$ is said the {\it degree of the curve} ${\mathcal
C}_l$. For a Riemann surface $S$, a well-known result is
($[2]$){\it there is a holomorphic mapping $\varphi : S\rightarrow
P_2{\mathcal C}$ such that $\varphi (S)$ is a complex plane
algebraic curve and}

$$g(S)=\frac{(d(\varphi (S))-1)(d(\varphi (S))-2)}{2}.$$

By map theory, we know a combinatorial map also is on a surface
with genus. Then {\it whether we can get an algebraic curve by all
edges in a map or by make operations on the vertices or edges of
the map to get plane algebraic curve with given $k$-multiple
points?} and {\it how do we find the equation $f(x,y,z)=0$?}

\vskip 4mm

{\bf $5.3$ Classification of $s$-manifolds by maps}

\vskip 3mm

We present an elementary classification for the closed
$s$-manifolds in the Section $3$. For the general $s$-manifolds,
their correspondence combinatorial model is the maps on surfaces
with boundary, founded by Bryant and Singerman in $1985$ ([1]).
The later is also related to the modular groups of spaces and need
to investigate further itself. The questions are\vskip 3mm

$(i)$ {\it how can we combinatorially classify the general
$s$-manifolds by maps with boundary?}

$(ii)$ {\it how can we find the automorphism group of an
$s$-manifold?}

$(iii)$ {\it how can we know the numbers of non-isomorphic
$s$-manifolds, with or without root?}

$(iv)$ {\it find rulers for drawing an $s$-manifold on a surface,
such as, the torus, the projective plane or Klein bottle, not the
plane.}\vskip 2mm

The $s$-manifolds only using the triangulations of surfaces with
vertex valency in $\{5,6,7\}$. Then {\it what are the geometrical
mean of the other maps, such as, the $4$-regular maps on
surfaces.} It is already known that the later is related to the
Gauss cross problem of curves($[9]$).

\vskip 4mm

{\bf $5.4$ Map geometries}

\vskip 3mm

As we have seen in the previous section, map geometries are the
nice model of the Smarandache geometries. More works should be
dong for them.

\vskip 3mm

($i$) {\it For a given graph, determine properties of the map
geometries underlying this graph.}

($ii$) {\it For a given locally orientable surface, determine the
properties of map geometries on this surface.}

($iii$) {\it Classify the map geometries on a locally orientable
surface.}

($iv$) {\it Enumerate non-equivalent map geometries underlying a
graph or on a locally orientable surface.}

($v$) {\it Establish the surface geometry by map geometries.}

\vskip 2mm

\vskip 4mm

{\bf $5.5$ Gauss mapping among surfaces}

\vskip 3mm

In the classical differential geometry, a {\it Gauss mapping}
among surfaces is defined as follows([10]):\vskip 3mm

{\it Let ${\mathcal S}\subset R^3$ be a surface with an
orientation {$\bf N$}. The mapping $ N: {\mathcal S}\rightarrow
R^3$ takes its value in the unit sphere}

$$S^2=\{(x,y,z)\in R^3|x^2+y^2+z^2=1\}$$

\no{\it along the orientation {$\bf N$}. The map $ N: {\mathcal
S}\rightarrow S^2$, thus defined, is called the Gauss
mapping.}\vskip 3mm

we know that for a point $P\in {\mathcal S}$ such that the
Gaussian curvature $K(P)\not=0$ and $V$ a connected neighborhood
of $P$ with $K$ does not change sign,

$$K(P)=\lim_{A\rightarrow 0}\frac{N(A)}{A},$$

\no{where $A$ is the area of a region $B\subset V$ and $ N(A)$ is
the area of the image of $B$ by the Gauss mapping $ N: {\mathcal
S}\rightarrow S^2$. The questions are}

{\it $(i)$ what is its combinatorial meaning of the Gauss mapping?
How to realizes it by maps?

$(ii)$ how we can define various curvatures for maps and rebuilt
the results in the classical differential geometry?}

\vskip 4mm

{\bf $5.6$ The Gauss-Bonnet theorem}

\vskip 3mm

{\it Let ${\mathcal S}$ be a compact orientable surface. Then}

$$\int\int_{\mathcal S}Kd\sigma =2\pi\chi ({\mathcal S}),$$

\no{\it where $K$ is Gaussian curvature on ${\mathcal S}$.}

 This is the famous {\it Gauss-Bonnet theorem} for compact surface
($[2],[6])$. The questions are

{\it $(i)$ what is its combinatorial mean of the Gauss curvature?}

{\it $(ii)$  how can we define the angle, area, volume, curvature,
$\cdots$,  of a map? }

{\it ($iii$) can we rebuilt the Gauss-Bonnet theorem by maps? or
can we get a generalization of the classical Gauss-Bonnet theorem
by maps? }

\vskip 4mm

{\bf $5.7$ Riemann manifolds}

\vskip 3mm

A Riemann surface is just a Riemann $2$-manifold, which has become
a source of the mathematical creative power. A {\it Riemann
$n$-manifold} $(M,g)$ is a $n$-manifold $M$ with a Riemann metric
$g$. Many important results in Riemann surfaces are generalized to
Riemann manifolds with higher dimension ($[6]$). For example, let
${\mathcal M}$ be a complete, simple-connected Riemann
$n$-manifold with constant sectional curvature $c$, then we know
that {\it ${\mathcal M}$ is isometric to one of the model spaces
${\mathcal R}^n, S_{{\mathcal R}^n}$ or $H_{{\mathcal R}^n}$}.
{\it Whether can we systematically rebuilt the Riemann manifold
theory by combinatorial maps? or  can we make a combinatorial
generalization of results in the Riemann geometry, for example,
the Chern-Gauss-Bonnet theorem ($[6]$)? }

\vskip 8mm

{\bf References}\vskip 5mm

\re{[1]}R.P.Bryant and D.Singerman, Foundations of the theory of
maps on surfaces with boundary,{\it
Quart.J.Math.Oxford}(2),36(1985), 17-41.

\re{[2]}H. M.Farkas and I. Kra, {\it Riemann Surfaces},
Springer-Verlag New York inc(1980).

\re{[3]}H.Iseri, {\it Smarandache manifolds}, American Research
Press, Rehoboth, NM,2002.

\re{[4]}H.Iseri, {\it Partially Paradoxist Smarandache
Geometries}, http://www.gallup.unm.
edu/\~smarandache/Howard-Iseri-paper.htm.

\re{[5]}L.Kuciuk and M.Antholy, An Introduction to Smarandache
Geometries, {\it Mathematics Magazine, Aurora, Canada},
Vol.12(2003)

\re{[6]}J.M.Lee, {\it Riemann Manifolds}, Springer-Verlag New
York,Inc(1997).

\re{[7]}Y.P.Liu, {\it Advances in Combinatorial Maps}, Northern
Jiaotong University Publisher, Beijing (2003).

\re{[8]}Yanpei Liu, {\it Enumerative Theory of Maps}, Kluwer
Academic Publisher, Dordrecht / Boston / London (1999).

\re{[9]}Y.P.Liu, {\it Embeddability in Graphs}, Kluwer Academic
Publisher, Dordrecht / Boston / London (1995).

\re{[10]}Mantredo P.de Carmao, {\it Differential Geometry of
Curves and Surfaces}, Pearson Education asia Ltd (2004).

\re{[11]}L.F.Mao, {\it Automorphism groups of maps, surfaces and
Smarandache geometries}, American Research Press, Rehoboth,
NM,2005.

\re{[12]}L.F.Mao and Y.P.Liu, A new approach for enumerating maps
on orientable surfaces, {\it Australasian J. Combinatorics},
vol.30(2004), 247-259.

\re{[13]}L.F.Mao and Y.P.Liu, Group action approach for
enumerating maps on surfaces,{\it J.Applied Math. \& Computing},
vol.13, No.1-2,201-215.

\re{[14]}W.S.Massey, {\it Algebraic topology: an introduction},
Springer-Verlag,New York, etc.(1977).

\re{[15]}V.V.Nikulin and I.R.Shafarevlch, {\it Geometries and
Groups}, Springer-Verlag Berlin Heidelberg (1987)

\re{[16]}F. Smarandache, Mixed noneuclidean geometries, {\it
eprint arXiv: math/0010119}, 10/2000.

\re{[17]}J.Stillwell, {\it Classical topology and combinatorial
group theory}, Springer-Verlag New York Inc., (1980).

\re{[18]} W.T.Tutte, What is a maps? in {\it New Directions in the
Theory of Graphs} (ed.by F.Harary), Academic Press (1973),
309~325.

\end{document}